\documentclass[11pt]{amsart}
\usepackage{amsfonts}
\usepackage{mathrsfs}
\usepackage{amssymb,latexsym}
\usepackage{pstcol,pstricks,color}

 \numberwithin{equation}{section}

\newtheorem{theorem}{Theorem}[section]

\newtheorem{corollary}[theorem]{Corollary}

\newtheorem{remark}[theorem]{Remark}

\newtheorem{lemma}[theorem]{Lemma}

\newtheorem{example}[theorem]{Example}

\def\qed{\hfill $\Box$}
\def\pf{\noindent {\it Proof.} }

\title{$\lambda$-factorials of $n$}

\begin{document}
\maketitle
\begin{center}
Yidong Sun$^\dag$ and Jujuan Zhuang$^\ddag$

%%%%\thanks{Corresponding author: Yidong Sun, sydmath@yahoo.com.cn.}%%

Department of Mathematics, Dalian Maritime University, 116026 Dalian, P.R. China\\[5pt]

{\it  Email: $^\dag$sydmath@yahoo.com.cn, $^\ddag$jjzhuang1979@yahoo.com.cn}

\end{center}\vskip0.5cm

\subsection*{Abstract}
Recently, by  the Riordan's identity related to tree enumerations,
\begin{eqnarray*}
\sum_{k=0}^{n}\binom{n}{k}(k+1)!(n+1)^{n-k}  &=& (n+1)^{n+1},
\end{eqnarray*}
Sun and Xu derived another analogous one,
\begin{eqnarray*}
\sum_{k=0}^{n}\binom{n}{k}D_{k+1}(n+1)^{n-k}  &=& n^{n+1},
\end{eqnarray*}
where $D_{k}$ is the number of permutations with no fixed points on $\{1,2,\dots, k\}$. In the
paper, we utilize the $\lambda$-factorials of $n$, defined by Eriksen, Freij and W$\ddot{a}$stlund, to give
a unified generalization of these two identities.
We provide for it a combinatorial proof by the functional digraph theory and another two algebraic proofs.
Using the umbral representation of our generalized identity and the Abel's binomial formula, we deduce several properties
for $\lambda$-factorials of $n$ and establish the curious relations between the generating functions of general and exponential
types for any sequence of numbers or polynomials.

\medskip

{\bf Keywords}: Derangement; $\lambda$-factorial of $n$; Charlier polynomial;
Bell polynomial; Hermite polynomial.

\noindent {\sc 2000 Mathematics Subject Classification}: Primary
05A05; Secondary 05C30

\section{Introduction}

Let $\mathfrak{S}_n$ denote the set of permutations of
$[n]=\{1,2,\dots, n\}$. A {\em fixed point} of a permutation $\pi\in
\mathfrak{S}_n$ is an element $i\in [n]$ such that $\pi(i) = i$.
Denote by $fix(\pi)$ the number of fixed points of $\pi$. Recently,
Eriksen, Freij and W$\ddot{a}$stlund \cite{Eriksen} defined the
polynomials, called the {\it $\lambda$-factorials of $n$}, namely
\begin{eqnarray*}
f_n(\lambda) = \sum_{\pi\in \mathfrak{S}_n}\lambda^{fix(\pi)}, \hskip.5cm f_0(\lambda)=1.
\end{eqnarray*}
They utilized them to give closed formulas for the number of derangements (permutations with no fixed points)
with descents in prescribed positions and derived several nice properties for $f_n(\lambda)$ such as
\begin{eqnarray}
f_n(\lambda+\mu) &=& \sum_{k=0}^{n}\binom{n}{k}f_k(\lambda)\mu^{n-k}, \label{eqn 1.0a}\\
f_n(\lambda)     &=& \sum_{k=0}^{n}\binom{n}{k}k!(\lambda-1)^{n-k},  \label{eqn 1.0b} \\
f_n(\lambda)     &=&  nf_{n-1}(\lambda)+(\lambda-1)^{n},              \label{eqn 1.0c} \\
\frac{d}{d\lambda}f_n(\lambda)     &=&  nf_{n-1}(\lambda).           \label{eqn 1.0d}
\end{eqnarray}

Clearly, we have $f_n(0) = D_n$ \cite[A000166]{Sloane} and $f_n(1) = n!$, where $D_n$ is the number of derangements
in $\mathfrak{S}_n$. The relation (\ref{eqn 1.0d}) indicates that $f_n(\lambda)$ $(n=0,1,\dots)$ form a kind of
special Appell polynomials \cite{Appell}. According to the definition, $f_n(\lambda)$ also has another expression
\begin{eqnarray}\label{eqn 1.0e}
f_n(\lambda) &=& \sum_{k=0}^{n}\binom{n}{k}D_{k}\lambda^{n-k}.
\end{eqnarray}
It should be noted that $f_n(\lambda)$ has close relation to the (re-normalized) Charlier
polynomials $C_n(\alpha, u)$ \cite{Gesselb} defined by
\begin{eqnarray*}
C_n(\alpha, u)=\sum_{k=0}^{n}\binom{n}{k}(\alpha)_{k}u^{n-k},
\end{eqnarray*}
where $(\alpha)_k=\alpha(\alpha+1)\cdots (\alpha+k-1)$. Clearly, $f_n(\lambda)=C_n(1,\lambda-1)$.

Recall that the Riordan identity in \cite[P173]{Comtet} and \cite{Riordan} states
\begin{eqnarray*}
\sum_{k=0}^{n}\binom{n}{k}(k+1)!(n+1)^{n-k}  &=& (n+1)^{n+1}.
\end{eqnarray*}
Using the Riordan identity, Sun and Xu \cite{SunXu} deduced an analogous identity,
\begin{eqnarray*}
\sum_{k=0}^{n}\binom{n}{k}D_{k+1}(n+1)^{n-k}  &=& n^{n+1}.
\end{eqnarray*}
Motivated by these two remarkable identities, we give the following general one and provide with a combinatorial
interpretation by the functional digraph theory.
\begin{theorem} For any integer $n\geq 0$ and any indeterminant $\lambda$, there holds
\begin{eqnarray}\label{eqn 1.1}
\sum_{k=0}^{n}\binom{n}{k}f_{k+1}(\lambda)(n+1)^{n-k} &=& (n+\lambda)^{n+1}.
\end{eqnarray}
\end{theorem}

Using the umbral representation of (\ref{eqn 1.1}) and the Abel's binomial formula, we have the second main result.
\begin{theorem} For any sequence $(a_n)_{n\geq 0}$, let $A(x)=\sum_{n\geq 0}a_n\frac{x^n}{n!}$.
Then
\begin{eqnarray}\label{eqn 1.2}
\sum_{n\geq 0}a_nf_n(\lambda)\frac{x^n}{n!} &=& \sum_{k\geq 0}(k+\lambda-1)^{k}\frac{x^kA^{(k)}(-kx)}{k!},
\end{eqnarray}
where $A^{(k)}(-kx)$ denotes the $k$-th derivative of $A(x)$ taking
value at $-kx$. In particular, the case $\lambda=1$ generates
\begin{eqnarray*}
\sum_{n\geq 0}a_nx^n  &=& \sum_{k\geq 0}\frac{(kx)^{k}A^{(k)}(-kx)}{k!}.
\end{eqnarray*}
\end{theorem}

The organization of this paper is as follows. The next section is
devoted to the proofs of (\ref{eqn 1.1}). In the third section,
we focus on the proof of (\ref{eqn 1.2}) and give some applications.
In the forth section, using the umbral representation of (\ref{eqn 1.1})
we further to consider the properties of $f_n(\lambda)$, involving
its another explicit formula. In the final section, we give some comments and
provide several open problems.

\section{Three Proofs for (\ref{eqn 1.1})}

In this section, we give three different proofs for (\ref{eqn 1.1}), one is a combinatorial proof by the functional digraph theory,
one is a generating function proof, and the other is a proof by umbral calculus.

\subsection{  First proof of (\ref{eqn 1.1}). } In order to give the combinatorial proof of (\ref{eqn 1.1}), we need some notations.
A {\it rooted labeled tree} on $[n]$ is an acyclic connected graph on the vertex set $[n]$ such
that one vertex, called the {\it root}, is distinguished. A {\it labeled forest} is a graph such that
every connected component is a rooted labeled tree. We denote by $\mathfrak{F}_n$ the set of labeled forests on $[n+1]$ and by
$\mathfrak{F}_{n,k}$ the set of labeled forests $[n+1]$ with exactly $k+1$ components. It is well known that
Cayley's formuals \cite{Cayley} state that $|\mathfrak{F}_n|=(n+2)^n$ and $|\mathfrak{F}_{n,k}|=\binom{n}{k}(n+1)^{n-k}$.

Let $\mathfrak{M}_n$ denote the set of maps $\sigma : [n]\rightarrow [n]$. Clearly $|\mathfrak{M}_n|=n^{n}$. For any $\sigma \in \mathfrak{M}_n$, we
represent $\sigma$ as a directed graph $G_{\sigma}$ by drawing arrows from $i$ to $\sigma(i)$. For any component of $G_{\sigma}$, it
contains equally many vertices and edges, and hence has exactly one directed cycle. Let $\mathfrak{R}$ denote the set of all the vertices of
these cycles of $G_{\sigma}$. Precisely, $\sigma|_{\mathfrak{R}}$, the restriction of $\sigma$ onto $\mathfrak{R}$ is just a permutation
on $\mathfrak{R}$. If deleting all the directed edges in these cycles, the remainders (omitting the directions, since all edges
are directed towards the roots) form a labeled forest $F$ on $[n]$. Conversely, it is not difficult to recover the map $\sigma$ from the pair
$(F,\pi)$, where $\pi$ is a permutation of the set $\mathfrak{R}_F$ of the roots of $F$.
Hence there exists a bijection between the set $\mathfrak{P}_n$ of the pairs $(F,\pi)$
and $\mathfrak{M}_{n+1}$, where $F\in \mathfrak{F}_n$. See \cite{Aigner, Joyal} for more details.

Now we can give a combinatorial interpretation for (\ref{eqn 1.1}).

It suffices to prove (\ref{eqn 1.1}) for the cases when $\lambda$ are nonnegative integers.
Let $\mathfrak{M}_{n+\lambda+1}^{*}$ be the set $\sigma\in \mathfrak{M}_{n+\lambda+1}$ such that $\sigma^{-1}(n+1)=\emptyset$ and
$\sigma(k)=k$ for $n+2\leq k\leq n+\lambda+1$. Clearly, $|\mathfrak{M}_{n+\lambda+1}^{*}|=(n+\lambda)^{n+1}$.

For any $\sigma\in \mathfrak{M}_{n+\lambda+1}^{*}$, it can uniquely determine a map $\tau$ from $[n+1]$ to $[n+1]$ such that the fixed points of $\tau$
have $\lambda$ colors, say, $c_1, c_2,\dots, c_{\lambda}$. The map $\tau$ is defined as follows.
\begin{eqnarray*}
\tau(i)=\left\{
\begin{array}{cll}
\sigma(i),     & if\ \sigma(i)\neq i\ and\ \sigma(i)\in [n],   \\
n+1  ,         & if\ \sigma(i)= i\ and\ i\in [n],      \\
i_{c_j},       & if\ \sigma(i)= n+j+1\ and\ i\in [n+1], j\in [\lambda],
\end{array}\right.
\end{eqnarray*}
where $\tau(i)=i_{c_j}$ means $\tau(i)=i$ and $i$ has color $c_{j}$. Conversely, one can uniquely recover $\sigma$ from $\tau$ by the following manner,
\begin{eqnarray*}
\sigma(i)=\left\{
\begin{array}{cll}
\tau(i),       & if\ \tau(i)\neq i\ and\ \tau(i)\in [n], \\
i,             & if\ \tau(i)= n+1 \ and\ i\in [n], \\
n+j+1 ,        & if\ \tau(i)=i_{c_j} \ and\ i\in [n+1], j\in [\lambda].
\end{array}\right.
\end{eqnarray*}

In other words, $G_{\tau}$ is obtained from $G_{\sigma}$ by the three steps:
\begin{itemize}
\item[(i)] Each directed cycle from $i$ to itself for $i\in [n]$ is transferred to be a directed edge from $i$ to $n+1$;
\item[(ii)] Each directed edge from $i$ to $n+j+1$ for $j\in [\lambda]$ is transferred to be a directed cycle from $i$ to itself for $i\in [n+1]$,
and such $i$ is assigned a color $c_j$;
\item[(iii)] Remove all the vertices $n+j+1$ for $j\in [\lambda]$.
\end{itemize}

It is clear that the procedure above is invertible and it is easy to recover $G_{\sigma}$ from $G_{\tau}$. So such maps $\tau$ are counted by
$(n+\lambda)^{n+1}$. On the other hand, we have also known that $\tau$ is bijected
to a pair $(F,\pi)\in \mathfrak{P}_n$ such that the fixed points of $\pi$ (also the fixed points of $\tau$) have $\lambda$ colors. If we restrict
$F\in \mathfrak{F}_{n,k}$, then $\pi$ is a permutation on $\mathfrak{R}_F$ with $k+1$ vertices such that the fixed points of $\pi$
have $\lambda$ colors. So such $F$ are counted by $|\mathfrak{F}_{n,k}|=\binom{n}{k}(n+1)^{n-k}$ and such $\pi$ are counted by $f_{k+1}(\lambda)$.
Summering all possible cases for $0\leq k\leq n$, we get (\ref{eqn 1.1}). \qed

\subsection{ Second proof of (\ref{eqn 1.1})} Let $y:=y(x)$ denote the exponential generating function for the labeled rooted trees
on $[n]$ which are counted by the sequence $(n^{n-1})_{n\geq 1}$, that is
$$y=\sum_{n\geq 1}n^{n-1}\frac{x^n}{n!}.$$
This generating function satisfies the relation $y=xe^{y}$ \cite{Stanley2}. By the Lagrange inversion formula, one can derive
\begin{eqnarray}
\frac{y^{k}}{k!} &=& \sum_{n\geq k}\binom{n-1}{k-1}n^{n-k}\frac{x^n}{n!},  \label{eqn 2.1} \\
\frac{e^{\lambda y}}{1-y} &=&  \sum_{n\geq 0}(n+\lambda)^{n}\frac{x^n}{n!}.   \label{eqn 2.2}
\end{eqnarray}

By (\ref{eqn 1.0b}), the exponential generating function $f(\lambda, t)$ for $f_n(\lambda)$ can be easily deduced
\begin{eqnarray}\label{eqn 2.3}
f(\lambda, t) &=& \sum_{k\geq 0}f_k(\lambda)\frac{t^k}{k!}=\frac{e^{(\lambda-1)t}}{1-t}.
\end{eqnarray}

Setting $t:=y$ in (\ref{eqn 2.3}), by (\ref{eqn 2.1}) and (\ref{eqn 2.2}), attracting the coefficient of $\frac{x^n}{n!}$, we have
\begin{eqnarray}\label{eqn 2.3a}
\sum_{k=1}^{n}\binom{n-1}{k-1}f_{k}(\lambda)n^{n-k} &=& (n+\lambda-1)^{n},
\end{eqnarray}
which is equivalent to (\ref{eqn 1.1}). \qed

\subsection{Third proof of (\ref{eqn 1.1})} Let $\mathbf{D}$ denote the umbra, given by $\mathbf{D}^n=D_n$.
See \cite{Gesselb, Roman, RomRota} for more information on umbral calculus. By (\ref{eqn 1.0e}),
$f_n(\lambda)$ can be represented umbrally as
\begin{eqnarray}\label{eqn 2.4}
f_{n}(\lambda) &=& (\mathbf{D}+\lambda)^{n}.
\end{eqnarray}
Then we have
\begin{eqnarray*}
\lefteqn {\sum_{k=0}^{n}\binom{n}{k}f_{k+1}(\lambda)(n+1)^{n-k}}   \\
&=& (\mathbf{D}+\lambda)(\mathbf{D}+\lambda+n+1)^{n} \\
&=& (\mathbf{D}+\lambda+n+1)^{n+1}-(n+1)(\mathbf{D}+\lambda+n+1)^{n}      \\
&=& f_{n+1}(\lambda+n+1)-(n+1)f_{n}( \lambda+n+1 )            \\
&=& (n+\lambda)^{n+1},               \hskip2cm (\mbox{by}\ (\ref{eqn 1.0c}))
\end{eqnarray*}
as desired. \qed

\section{Proof of Theorem 1.2 and Its Applications}

In this section, we first give the proof of Theorem 1.2, and then we provide several interesting examples.

\subsection{Proof of Theorem 1.2} Recall the Abel binomial theorem \cite[P128]{Comtet} states
\begin{eqnarray}\label{eqn 3.1}
 (a+b)^n &=& \sum_{k=0}^{n}\binom{n}{k}a(a-kt)^{k-1}(b+kt)^{n-k}.
\end{eqnarray}

Let $A(x)=\sum_{n\geq 0}a_n\frac{x^n}{n!}$ be the exponential generating function for any sequence $(a_n)_{n\geq 0}$, then the
Abel binomial theorem is equivalent to the form \cite[P130]{Comtet}
\begin{eqnarray}\label{eqn 3.2}
A(x) &=& \sum_{k\geq 0}x(x-kt)^{k-1}\frac{A^{(k)}(kt)}{k!},
\end{eqnarray}
where $t$ is a new indeterminant and $A^{(k)}(kt)$ denotes the $k$-th derivative of $A(x)$ taking
value at $kt$.

By (\ref{eqn 1.1}) and (\ref{eqn 2.4}), we have
\begin{eqnarray}\label{eqn 3.3}
(\mathbf{D}+\lambda)(\mathbf{D}+\lambda+n+1)^{n} &=& (n+\lambda)^{n+1}.
\end{eqnarray}

Setting $x:=(\mathbf{D}+\lambda)x$ and $t=-x$ in (\ref{eqn 3.2}), by (\ref{eqn 3.3}), one can obtain (\ref{eqn 1.2}).

\subsection{Applications of Theorem 1.2}

In this subsection, as applications of Theorem 1.2, we only consider three special cases when $a_n$ are taken to be the Charlier, Bell
and Hermite polynomials. Of course, one can also consider other interesting cases such as $a_n$ are the Bessel,
Chebyshev, Legendre, Jacobi, Laguerre, and ultraspherical polynomials and so on.

\begin{example}
Let $a_n=C_n(\alpha,u)$, the (re-normalized) Charlier polynomial, which has the exponential generating function $A(x)=\frac{e^{ux}}{(1-x)^{\alpha}}$.
It is easy to derive that $\frac{\partial}{\partial x}A(x)=\frac{\alpha+u(1-x)}{1-x}A(x)$, which implies the recurrence relation
\begin{eqnarray*}
C_{n+1}(\alpha,u) &=& \alpha C_{n}(\alpha+1,u)+uC_{n}(\alpha,u).
\end{eqnarray*}
Using this recurrence and by induction on $k$, one can deduce that
\begin{eqnarray*}
\frac{\partial^{k}}{\partial x^k}A(x)=\frac{C_{k}(\alpha,u(1-x))}{(1-x)^{k}}A(x)=\frac{C_{k}(\alpha,u(1-x))e^{ux}}{(1-x)^{\alpha+k}}.
\end{eqnarray*}
Then by Theorem 1.1, we have
\begin{eqnarray}\label{eqn 3.4}
\sum_{n\geq 0}C_n(\alpha,u)f_n(\lambda)\frac{x^n}{n!} &=& \sum_{k\geq 0}\frac{(k+\lambda-1)^{k}x^kC_{k}(\alpha,u(1+kx))e^{-ukx}}{k!(1+kx)^{\alpha+k}}.
\end{eqnarray}

More generally, if let $a_n=C_{m+n}(\alpha,u)$ or $A(x)=\frac{\partial^m}{\partial x^m}\frac{e^{ux}}{(1-x)^{\alpha}}$, then
\begin{eqnarray*}
\frac{\partial^{k}}{\partial x^k}A(x)=\frac{\partial^{m+k}}{\partial x^{m+k}}\frac{e^{ux}}{(1-x)^{\alpha}}=\frac{C_{m+k}(\alpha,u(1-x))e^{ux}}{(1-x)^{\alpha+m+k}}.
\end{eqnarray*}
In this case Theorem 1.1 generates that
\begin{eqnarray}\label{eqn 3.5}
\sum_{n\geq 0}C_{m+n}(\alpha,u)f_n(\lambda)\frac{x^n}{n!} &=& \sum_{k\geq 0}\frac{(k+\lambda-1)^{k}x^kC_{m+k}(\alpha,u(1+kx))e^{-ukx}}{k!(1+kx)^{\alpha+m+k}}.
\end{eqnarray}

The parameter specializations in (\ref{eqn 3.4}) and (\ref{eqn 3.5}) produce several consequences.

\begin{itemize}
\item[Case 1.] When $\alpha=1, u=\mu-1$, $C_{m+n}(1, \mu-1)=f_{m+n}(\mu)$. Then by (\ref{eqn 3.5}) we have
\begin{eqnarray*}
\sum_{n\geq 0}f_{m+n}(\mu)f_n(\lambda)\frac{x^n}{n!} &=& \sum_{k\geq 0}\frac{(k+\lambda-1)^{k}x^kf_{m+k}(1+(\mu-1)(1+kx))e^{-(\mu-1)kx}}{k!(1+kx)^{m+k+1}}.
\end{eqnarray*}
which, when $\mu=\lambda=0$, by $f_{n}(0)=D_{n}$, yields
\begin{eqnarray}\label{eqn 3.6}
\sum_{n\geq 0}D_{m+n}D_n \frac{x^n}{n!}  &=& \sum_{k\geq 0} \frac{(k-1)^{k}x^kf_{m+k}(-kx)e^{kx}}{k!(1+kx)^{m+k+1}},
\end{eqnarray}
and when $\lambda=1$ leads to
\begin{eqnarray}\label{eqn 3.7}
\sum_{n\geq 0}f_{m+n}(\mu)x^n  &=& \sum_{k\geq 0}\frac{(kx)^{k}f_{m+k}(1+(\mu-1)(1+kx))e^{-(\mu-1)kx}}{k!(1+kx)^{m+k+1}}.
\end{eqnarray}
The case $\mu=0$ in (\ref{eqn 3.7}) gives the ordinary generating function for $D_{m+n}$,
\begin{eqnarray*}
\sum_{n\geq 0}D_{m+n}x^n  &=& \sum_{k\geq 0}\frac{(kx)^{k}f_{m+k}(-kx)e^{kx}}{k!(1+kx)^{m+k+1}}.
\end{eqnarray*}

\item[Case 2.] When $u=0$, $C_{m+n}(\alpha,0)=(\alpha)_{m+n}$. Then by (\ref{eqn 3.5}) we have
\begin{eqnarray*}
\sum_{n\geq 0}(\alpha)_{m+n}f_n(\lambda)\frac{x^n}{n!} &=& \sum_{k\geq 0}\frac{(\alpha)_{m+k}(k+\lambda-1)^{k}x^k}{k!(1+kx)^{m+k+1}},
\end{eqnarray*}
which, when $\alpha=1, \lambda=1$ and $\alpha=1, m=0$, leads respectively to the ordinary generating function for $(m+n)!$ and $f_n(\lambda)$,
\begin{eqnarray}\label{eqn 3.7.1}
\sum_{n\geq 0}(m+n)!x^n        &=&  \sum_{k\geq 0}\frac{(m+k)!}{k!}\frac{(kx)^{k}}{(1+kx)^{m+k+1}}, \nonumber\\
\sum_{n\geq 0}f_n(\lambda)x^n  &=&  \sum_{k\geq 0}\frac{(k+\lambda-1)^{k}x^k}{(1+kx)^{k+1}}.
\end{eqnarray}

\end{itemize}

\end{example}

\begin{remark}
Gessel \cite{Gesselb} utilized the umbral calculus method to derive the bilinear generating function for Charlier polynomials
\begin{eqnarray*}
\sum_{n\geq 0}C_n(\alpha,u)C_n(\beta,v)\frac{x^n}{n!}  &=& e^{uvx}\sum_{k\geq 0}
\frac{(\alpha)_k}{(1-vx)^{\alpha+k}}\frac{(\beta)_k}{(1-ux)^{\beta+k}}\frac{x^{k}}{k!},
\end{eqnarray*}
which, when $\alpha=1, u=-1$, gives us another more interesting but considerably more recondite formula analogous to the case $m=0$ in (\ref{eqn 3.6}),
\begin{eqnarray*}
\sum_{n\geq 0}D_{n}D_n \frac{x^n}{n!}  &=& e^{x}\sum_{k\geq 0} \frac{k!x^{k}}{(1+x)^{2k+2}}.
\end{eqnarray*}
\end{remark}

\begin{remark}
Clarke, Han and Zeng \cite{ClarkHanZeng} utilized the Laplace transformation to deduce another ordinary generating function for $f_{n}(\mu)$ analogous
to (\ref{eqn 3.7.1}) or the case $m=0$ in (\ref{eqn 3.7}),
\begin{eqnarray*}
\sum_{n\geq 0}f_{n}(\mu)x^n  &=& \sum_{k\geq 0}\frac{k!x^k}{(1-(\mu-1)x)^{k+1}}.
\end{eqnarray*}
\end{remark}

\begin{example}
Let $a_n=B_{n}(u)$, the $n$-the Bell polynomial, which has the exponential generating function $A(x)=\exp(u(e^{x}-1))$.
It is easy to derive that $\frac{\partial}{\partial x}A(x)=ue^xA(x)$, which implies the recurrence relation
\begin{eqnarray*}
B_{n+1}(u) &=& uB_{n}(u)+u\frac{d}{du}B_{n}(u).
\end{eqnarray*}
Using this recurrence and by induction on $k$, one can deduce that
\begin{eqnarray*}
\frac{\partial^{k}}{\partial x^k}A(x)= B_k(ue^x)A(x)=B_k(ue^x)\exp(u(e^{x}-1)).
\end{eqnarray*}
Then by Theorem 1.1, we have
\begin{eqnarray*}
\sum_{n\geq 0}B_n(u)f_n(\lambda)\frac{x^n}{n!} &=& \sum_{k\geq 0}\frac{(k+\lambda-1)^{k}x^kB_{k}(ue^{-kx})\exp{(u(e^{-kx}-1))}}{ k!}.
\end{eqnarray*}
More generally, if let $a_n=B_{m+n}(u)$ or $A(x)=\frac{\partial^m}{\partial x^m}\exp(u(e^{x}-1))$, then
\begin{eqnarray*}
\frac{\partial^{k}}{\partial x^k}A(x)=\frac{\partial^{m+k}}{\partial x^{m+k}}\exp(u(e^{x}-1))=B_{m+k}(ue^x)\exp(u(e^{x}-1)).
\end{eqnarray*}
In this case Theorem 1.1 generates that
\begin{eqnarray*}
\sum_{n\geq 0}B_{m+n}(u)f_n(\lambda)\frac{x^n}{n!} &=& \sum_{k\geq 0}\frac{(k+\lambda-1)^{k}x^kB_{m+k}(ue^{-kx})\exp{(u(e^{-kx}-1))}}{ k!},
\end{eqnarray*}
which, when $u=\lambda=1$, leads to the ordinary generating function for the Bell numbers $B_{m+n}(1)=B_{m+n}$ \cite[A000110]{Sloane},
\begin{eqnarray}\label{eqn 3.8}
\sum_{n\geq 0}B_{m+n}x^n   &=& \sum_{k\geq 0}\frac{(kx)^{k}B_{m+k}(e^{-kx})\exp{(e^{-kx}-1)}}{ k!}.
\end{eqnarray}

\end{example}

\begin{remark}
Another classical ordinary generating function for the Bell numbers $B_n$ is
\begin{eqnarray*}
\sum_{n\geq 0}B_{n}x^n &=& \sum_{k\geq 0}\frac{x^k}{(1-x)(1-2x)\cdots (1-kx)}.
\end{eqnarray*}
Klazar \cite{Klazar} in depth investigated this generating function and proved that it satisfies no
algebraic differential equation over the complex field.
\end{remark}

\begin{example}
Let $a_{n}=H_n(u)$, the $n$-th (re-normalized) Hermite polynomial, whose exponential generating function is
$A(x)=\exp(ux+\frac{x^2}{2})$. The polynomial $H_n(u)$ also counts involutions on $[n]$ such that the fixed points have $u$ colors.
It is easy to derive that $\frac{\partial}{\partial x}A(x)=(u+x)A(x)$, which implies the recurrence relation
\begin{eqnarray*}
H_{n+1}(u) &=& uH_{n}(u)+ \frac{d}{du}H_{n}(u).
\end{eqnarray*}
Using this recurrence and by induction on $k$, one can deduce that
\begin{eqnarray*}
\frac{\partial^{k}}{\partial x^k}A(x)= H_k(u+x)A(x)=H_k(u+x)\exp(ux+\frac{x^2}{2}).
\end{eqnarray*}
Then by Theorem 1.1, we have
\begin{eqnarray*}
\sum_{n\geq 0}H_n(u)f_n(\lambda)\frac{x^n}{n!} &=& \sum_{k\geq 0}\frac{(k+\lambda-1)^{k}x^kH_{k}(u-kx)\exp{(-ukx+\frac{(kx)^2}{2})}}{ k!}.
\end{eqnarray*}
More generally, if let $a_n=H_{m+n}(u)$ or $A(x)=\frac{\partial^m}{\partial x^m}\exp(ux+\frac{x^2}{2})$, then
\begin{eqnarray*}
\frac{\partial^{k}}{\partial x^k}A(x)=\frac{\partial^{m+k}}{\partial x^{m+k}}\exp(ux+\frac{x^2}{2})=H_{m+k}(u+x)\exp(ux+\frac{x^2}{2}).
\end{eqnarray*}
In this case Theorem 1.1 generates that
\begin{eqnarray}\label{eqn 3.9}
\sum_{n\geq 0}H_{m+n}(u)f_n(\lambda)\frac{x^n}{n!} &=& \sum_{k\geq 0}\frac{(k+\lambda-1)^{k}x^kH_{m+k}(u-kx)\exp{(-ukx+\frac{(kx)^2}{2})}}{ k!}.
\end{eqnarray}
The cases when $\lambda=1$ and $u=1$ or $u=0$ in (\ref{eqn 3.9}), lead respectively to the ordinary generating functions for
the involution numbers $I_{m+n}=H_{m+n}(1)$ \cite[A000085]{Sloane} and the matching numbers $M_{m+n}=H_{m+n}(0)$ \cite[A001147]{Sloane},
\begin{eqnarray*}
\sum_{n\geq 0}I_{m+n}x^n   &=& \sum_{k\geq 0}\frac{(kx)^{k} H_{m+k}(1-kx)\exp{(-kx+\frac{(kx)^2}{2})}}{ k!}, \\
\sum_{n\geq 0}M_{m+n}x^n   &=& \sum_{k\geq 0}\frac{(kx)^{k} H_{m+k}(-kx)\exp{(\frac{(kx)^2}{2})}}{ k!}.
\end{eqnarray*}

\end{example}

\section{Further properties of $f_n(\lambda)$}

\begin{theorem} For any integer $n\geq 0$ and any indeterminants $\lambda, \mu$, there hold
\begin{eqnarray}
\sum_{k=0}^{n}\binom{n}{k}f_{k}(\lambda)(\mu+k-n)\mu^{n-k} &=& \mu(\lambda+\mu-1)^{n},    \label{eqn 4.1} \\
\sum_{k=0}^{n}\binom{n}{k}f_{k}(\lambda)f_{n-k}(\mu+1) &=&  (\lambda+\mu-1)^{n+1}+(n-\lambda-\mu+2)f_n(\lambda+\mu).   \label{eqn 4.2}
\end{eqnarray}
\end{theorem}
\pf For (\ref{eqn 4.1}), we have
\begin{eqnarray*}
\lefteqn{ \sum_{k=0}^{n}\binom{n}{k}f_{k}(\lambda)(\mu+k-n)\mu^{n-k} } \\
&=& \mu f_{n}(\lambda+\mu)-\mu\frac{\partial }{\partial \mu} f_{n}(\lambda+\mu) \hskip1.1cm   (\mbox{by}\ (\ref{eqn 1.0a})) \\
&=& \mu f_{n}(\lambda+\mu)-n\mu  f_{n-1}(\lambda+\mu)   \hskip1cm   (\mbox{by}\ (\ref{eqn 1.0d}))  \\
&=& \mu(\lambda+\mu-1)^{n}.   \hskip3.45cm   (\mbox{by}\ (\ref{eqn 1.0c}))
\end{eqnarray*}
For (\ref{eqn 4.2}), setting $\mu:=\mathbf{D}+\mu+n+1$ in (\ref{eqn 4.1}), then the left hand side of (\ref{eqn 4.1}) equals
\begin{eqnarray}
LHS\ of\ (\ref{eqn 4.1}) &=& \sum_{k=0}^{n}\binom{n}{k}f_{k}(\lambda)( \mathbf{D}+\mu+n+1+k-n)( \mathbf{D}+\mu+n+1)^{n-k} \nonumber\\
&=& \sum_{k=0}^{n}\binom{n}{k}f_{k}(\lambda)\big(f_{n-k+1}(\mu+n+1)-(n-k)f_{n-k}(\mu+n+1)\big)   \hskip.4cm  (\mbox{by}\ (\ref{eqn 2.4})) \nonumber\\
&=& \sum_{k=0}^{n}\binom{n}{k}f_{k}(\lambda)\big(f_{n-k}(\mu+n+1)+(\mu+n)^{n-k+1}\big)  \hskip2cm   (\mbox{by}\ (\ref{eqn 1.0c}))  \nonumber\\
&=& \sum_{k=0}^{n}\binom{n}{k}f_{k}(\lambda)f_{n-k}(\mu+n+1)+(\mu+n)f_{n}(\lambda+\mu+n),    \hskip.8cm   (\mbox{by}\ (\ref{eqn 1.0a}))  \label{eqn 4.2a}
\end{eqnarray}
and the right hand side of (\ref{eqn 4.1}) is equal to
\begin{eqnarray}
RHS\ of\ (\ref{eqn 4.1}) &=& ( \mathbf{D}+\mu+n+1)( \mathbf{D}+\lambda+\mu+n)^{n} \nonumber \\
&=&  ( \mathbf{D}+\lambda+\mu-1)( \mathbf{D}+\lambda+\mu+n)^{n} + (n-\lambda+2)( \mathbf{D}+\lambda+\mu+n)^{n}   \nonumber \\
&=& (\lambda+\mu+n-1)^{n+1} +(n-\lambda+2)f_n(\lambda+\mu+n).   \hskip1cm   (\mbox{by}\ (\ref{eqn 1.0c})) \label{eqn 4.2b}
\end{eqnarray}
Then (\ref{eqn 4.2}) is followed by setting $\mu:=\mu-n$ in (\ref{eqn 4.2a}) and (\ref{eqn 4.2b}). \qed

\begin{remark}
It should be noted that (\ref{eqn 4.1}) in the case $\mu=1$ and (\ref{eqn 1.0e}) form a new inverse relation. In general, for any two
sequences $(a_n)_{n\geq 0}$ and $(b_n)_{n\geq 0}$,
\begin{eqnarray*}
b_n = \sum_{k=0}^{n}\binom{n}{k}D_{n-k}a_k &\Leftrightarrow & a_n = \sum_{k=0}^{n}\binom{n}{k}(1+k-n)b_{k}.
\end{eqnarray*}
\end{remark}

When $\lambda+\mu=n+2$, (\ref{eqn 4.2}) reduces to the surprising result.
\begin{corollary} For any integer $n\geq 0$ and any indeterminant $\lambda$, there holds
\begin{eqnarray*}
\sum_{k=0}^{n}\binom{n}{k}f_{k}(\lambda)f_{n-k}(n-\lambda+3) &=&  (n+1)^{n+1}.
\end{eqnarray*}
\end{corollary}

\begin{theorem} For any integer $n\geq 0$ and any indeterminants $\lambda, \mu$, there holds
\begin{eqnarray}\label{eqn 4.3}
 f_{n}(\lambda+\mu) &=& \sum_{k=0}^{n}\binom{n}{k}(\lambda+k)^{k}(\mu-k-1)^{n-k}.
\end{eqnarray}
\end{theorem}
\pf Setting $t=-1, a=\mathbf{D}+\lambda+1, b=\mu-1$ in (\ref{eqn 3.1}), we have
\begin{eqnarray*}
(\mathbf{D}+\lambda+\mu)^n &=& \sum_{k=0}^{n}\binom{n}{k}(\mathbf{D}+\lambda+1)(\mathbf{D}+\lambda+1+k)^{k-1}(\mu-k-1)^{n-k},
\end{eqnarray*}
which, by (\ref{eqn 2.4}), is equivalent to (\ref{eqn 4.3}). \qed
\begin{remark}
Note that (\ref{eqn 4.3}) in the case $\lambda:=\lambda-1, \mu=1$ and (\ref{eqn 2.3a}) form another known inverse
relation \cite[P164]{Comtet}. In general, for any two sequences $(a_n)_{n\geq 1}$ and $(b_n)_{n\geq 1}$,
\begin{eqnarray*}
b_n = \sum_{k=1}^{n}\binom{n-1}{k-1}n^{n-k}a_k &\Leftrightarrow & a_n = \sum_{k=1}^{n}(-1)^{n-k}\binom{n}{k}k^{n-k}b_{k}.
\end{eqnarray*}
\end{remark}

When $\mu=1-\lambda$ in (\ref{eqn 4.3}), by $f_n(1)=n!$, we have the well-known difference identity
\begin{eqnarray*}
\sum_{k=0}^{n}(-1)^{n-k}\binom{n}{k}(\lambda+k)^{n}  &=& n!.
\end{eqnarray*}

When $\mu=-\lambda$ in (\ref{eqn 4.3}), by $f_n(0)=D_n$, we have the following result.
\begin{corollary} For any integer $n\geq 0$ and any indeterminant $\lambda$, there holds
\begin{eqnarray} \label{eqn 4.3a}
D_n &=& \sum_{k=0}^{n}(-1)^{n-k}\binom{n}{k}(\lambda+k)^{k}(\lambda+k+1)^{n-k}.
\end{eqnarray}
\end{corollary}
\begin{remark}
The case $\lambda=-1$ in (\ref{eqn 4.3a}) has been obtained in \cite[P201]{Comtet} by the permanent theory. In fact, $D_n$ is also
the permanent of the matrix $\mathbf{J}-\mathbf{I}$, where $\mathbf{I}$ is the $n\times n$ unit matrix and $\mathbf{J}$ is the
$n\times n$ matrix with all entries being equal to $1$.
\end{remark}

\begin{theorem} For any integer $n\geq 0$ and any indeterminants $\lambda, \mu$, there hold
\begin{eqnarray}
\sum_{k=0}^{n}\binom{n}{k}f_{k+1}(\lambda)\mu^{n-k}      &=& \sum_{k=0}^{n}\binom{n}{k}(\lambda+k)^{k+1}(\mu-(n+1))(\mu-k-1)^{n-k-1},  \label{eqn 4.4} \\
\sum_{k=0}^{n}\binom{n}{k}f_{k+1}(\lambda)f_{n-k}(\mu+1) &=& \sum_{k=0}^{n}\binom{n}{k}(\lambda+k)^{k+1}(\mu-k-1)^{n-k}. \label{eqn 4.5}
\end{eqnarray}
\end{theorem}
\pf For (\ref{eqn 4.4}), setting $t=-1, a=\mu-n-1, b=\mathbf{D}+\lambda+n+1$ in (\ref{eqn 3.1}), we have
\begin{eqnarray*}
(\mathbf{D}+\lambda+\mu)^n &=& \sum_{k=0}^{n}\binom{n}{k} (\mathbf{D}+\lambda+n-k+1)^{n-k}(\mu-n-1)(\mu-n+k-1)^{k-1} \\
                           &=& \sum_{k=0}^{n}\binom{n}{k} (\mathbf{D}+\lambda+k+1)^{k}(\mu-n-1)(\mu-k-1)^{n-k-1}.
\end{eqnarray*}
Then, by (\ref{eqn 2.4}), we get
\begin{eqnarray*}
\lefteqn{ \sum_{k=0}^{n}\binom{n}{k}f_{k+1}(\lambda)\mu^{n-k} } \\
&=& (\mathbf{D}+\lambda)(\mathbf{D}+\lambda+\mu)^{n} \\
&=& \sum_{k=0}^{n}\binom{n}{k} (\mathbf{D}+\lambda)(\mathbf{D}+\lambda+k+1)^{k}(\mu-n-1)(\mu-k-1)^{n-k-1} \\
&=& \sum_{k=0}^{n}\binom{n}{k}(\lambda+k)^{k+1}(\mu-(n+1))(\mu-k-1)^{n-k-1}.
\end{eqnarray*}

For (\ref{eqn 4.5}), setting $\mu:=\mathbf{D}+\mu+n+1$ in (\ref{eqn 4.4}), by (\ref{eqn 2.4}), we have
\begin{eqnarray*}
\lefteqn { \sum_{k=0}^{n}\binom{n}{k}f_{k+1}(\lambda)f_{n-k}(\mu+n+1) } \\
&=& \sum_{k=0}^{n}\binom{n}{k}f_{k+1}(\lambda)(\mathbf{D}+\mu+n+1)^{n-k}   \\
&=& \sum_{k=0}^{n}\binom{n}{k}(\lambda+k)^{k+1}(\mathbf{D}+\mu)(\mathbf{D}+\mu+n-k)^{n-k-1} \\
&=& \sum_{k=0}^{n}\binom{n}{k}(\lambda+k)^{k+1}(\mu+n-k-1)^{n-k},
\end{eqnarray*}
which, by setting $\mu:=\mu-n$, generates (\ref{eqn 4.5}). \qed

\begin{remark}
The case $\mu=n+1$ in (\ref{eqn 4.4}) produces (\ref{eqn 1.1}).
\end{remark}

Setting $\mu=1-\lambda$ in (\ref{eqn 4.4}) and (\ref{eqn 4.5}), using the general difference identity \cite{Stanley1}
\begin{eqnarray*}
\sum_{k=0}^{n}(-1)^{n-k}\binom{n}{k}(\lambda+k)^{m}  &=& \sum_{k=n}^{m}(-)^{k}S(m,k)(-k)_{n}(-\lambda)_{k-n},
\end{eqnarray*}
where $S(m,k)$ is the Stirling number of the second kind \cite[A008277]{Sloane}, and by $S(n+1,n)=\binom{n+1}{2}$, we have
\begin{corollary} For any integer $n\geq 0$ and any indeterminant $\lambda$, there hold
\begin{eqnarray*}
\sum_{k=0}^{n}\binom{n}{k}f_{k+1}(\lambda)(1-\lambda)^{n-k}  &=& n!(\lambda+n),  \\
\sum_{k=0}^{n}\binom{n}{k}f_{k+1}(\lambda)f_{n-k}(2-\lambda) &=& (n+1)!(\lambda+\frac{n}{2}).
\end{eqnarray*}
\end{corollary}

\section{Comments and open questions}

In general, we can consider the generalization of (\ref{eqn 1.0a}) and (\ref{eqn 4.4}), that is
\begin{eqnarray*}
Q_{n,m}(\lambda, \mu)=\sum_{k=0}^{n}\binom{n}{k}f_{k+m}(\lambda)\mu^{n-k}.
\end{eqnarray*}
By (\ref{eqn 1.0c}), one can deduce the first recurrence relation for $Q_{n,m}(\lambda, \mu)$,
\begin{eqnarray}\label{eqn 5.1}
Q_{n,m}(\lambda, \mu)= nQ_{n-1,m}(\lambda, \mu)+mQ_{n,m-1}(\lambda, \mu)+(\lambda-1)^{m}(\lambda+\mu-1)^n
\end{eqnarray}
with the initial conditions $Q_{0,0}(\lambda, \mu)=1, Q_{n,0}(\lambda, \mu)=Q_{0,m}(\lambda, \mu)=0$ whenever $n,m<0$.
Clearly, (\ref{eqn 5.1}) reduces to (\ref{eqn 1.0c}) when $n=0$ and $m:=n$ or $m=0$ and $\mu=0$.

Let $Q(\lambda, \mu; x, t)$ denote the exponential generating function for $Q_{n,m}(\lambda, \mu)$, i.e.,
\begin{eqnarray*}
Q(\lambda, \mu; t, x) &=& \sum_{n,m\geq 0} Q_{n,m}(\lambda, \mu)\frac{t^n}{n!}\frac{x^m}{m!}.
\end{eqnarray*}
From (\ref{eqn 5.1}), we can derive the explicit formula for $Q(\lambda, \mu; t, x)$,
\begin{eqnarray}\label{eqn 5.2}
Q(\lambda, \mu; t, x) &=& \frac{e^{(\lambda+\mu-1)t}e^{(\lambda-1)x}}{1-t-x}.
\end{eqnarray}

By (\ref{eqn 5.2}), one has
\begin{eqnarray*}
\frac{\partial Q(\lambda, \mu; t, x)}{\partial t} &=& \frac{\partial Q(\lambda, \mu; t, x)}{\partial x}+\mu Q(\lambda, \mu; t, x),
\end{eqnarray*}
which implies that there holds another recurrence relation for $Q_{n,m}(\lambda, \mu)$,
\begin{eqnarray*}
Q_{n+1,m}(\lambda, \mu)= Q_{n,m+1}(\lambda, \mu)+\mu Q_{n,m}(\lambda, \mu).
\end{eqnarray*}

Note that the type of the exponential generating function $Q(\lambda, \mu; t, x)$ brings it into the general framework considered in \cite{SunXu},
which signifies that $Q_{n,m}(\lambda, \mu)$ has many other interesting properties. For examples, setting $t:=tx$ in (\ref{eqn 5.2}) and comparing
the coefficients of $\frac{x^N}{N!}$, we get
\begin{eqnarray*}
\sum_{n=0}^{N}\binom{N}{n}Q_{N-n,n}(\lambda, \mu)t^{N-n}= (t+1)^{N}f_{N}(\lambda+\frac{\mu t}{t+1}).
\end{eqnarray*}
Using the series expansion, we have
\begin{eqnarray*}
Q(\lambda, \mu; t, x) &=& \frac{e^{(\lambda+\mu-1)t}e^{(\lambda-1)x}}{1-t-x}=\frac{e^{(\lambda+\mu-1)t}}{1-t}\frac{e^{(\lambda-1)x}}{1-\frac{x}{1-t}} \\
                      &=& \sum_{m\geq 0}\frac{x^m}{m!}\sum_{j=0}^{m}\binom{m}{j}j!(\lambda-1)^{m-j}\frac{e^{(\lambda+\mu-1)t}}{(1-t)^{j+1}} \\
                      &=& \sum_{m\geq 0}\frac{x^m}{m!}\sum_{j=0}^{m}\binom{m}{j}j!(\lambda-1)^{m-j}\sum_{n\geq 0}\frac{t^n}{n!}\sum_{k=0}^{n}\binom{n}{k}(j+1)_k(\lambda+\mu-1)^{n-k} \\
                      &=& \sum_{n\geq 0}\sum_{m\geq 0}\frac{t^n}{n!}\frac{x^m}{m!}\sum_{k=0}^{n}\sum_{j=0}^{m}\binom{n}{k}\binom{m}{j}(k+j)!(\lambda+\mu-1)^{n-k}(\lambda-1)^{m-j}.
\end{eqnarray*}
Comparing the coefficients of $\frac{t^nx^m}{n!m!}$, we get an explicit formula for $Q_{n,m}(\lambda, \mu)$,
\begin{eqnarray*}
Q_{n,m}(\lambda, \mu) &=& \sum_{k=0}^{n}\sum_{j=0}^{m}\binom{n}{k}\binom{m}{j}(k+j)!(\lambda+\mu-1)^{n-k}(\lambda-1)^{m-j}.
\end{eqnarray*}

But here we have more interest in the type of formulas for $Q_{n,m}(\lambda, \mu)$ similar to (\ref{eqn 4.3}) and (\ref{eqn 4.4}).
\begin{lemma}
For any integer $n\geq 0$ and any indeterminants $\lambda, \mu$, there holds
\begin{eqnarray}\label{eqn 5.3}
Q_{n,m}(\lambda, \mu) &=& mQ_{n,m-1}(\lambda, \mathbf{D}+\mu+1)+(\lambda-1)^{m}f_n(\lambda+\mu),
\end{eqnarray}
or equivalently
\begin{eqnarray}\label{eqn 5.4}
Q_{n,m}(\lambda, \mu) &=& m\sum_{k=0}^{n}\binom{n}{k}f_{k+m-1}(\lambda)f_{n-k}(\mu+1)+(\lambda-1)^{m}f_n(\lambda+\mu).
\end{eqnarray}
\end{lemma}
\pf Note that
\begin{eqnarray*}
\lefteqn{\sum_{n\geq 0}Q_{n,m}(\lambda, \mu)\frac{t^n}{n!} = \frac{\partial^m}{\partial x^m} Q(\lambda, \mu; t, x)\Big|_{x=0}} \\
&=& \frac{\partial^m}{\partial x^m} \frac{e^{\mu t}e^{(\lambda-1)(t+x)}}{1-(t+x)}\Big|_{x=0}=e^{\mu t}\frac{\partial^m}{\partial t^m}\frac{e^{(\lambda-1)t}}{1-t}  \\
&=& f_m((\lambda-1)(1-t)+1)\frac{e^{(\lambda+\mu-1)t}}{(1-t)^{m+1}} \hskip3.5cm (\mbox{by\ Example\ 3.1})  \\
&=& mf_{m-1}((\lambda-1)(1-t)+1)\frac{e^{(\lambda+\mu-1)t}}{(1-t)^{m+1}}+(\lambda-1)^m\frac{e^{(\lambda+\mu-1)t}}{1-t} \hskip.5cm(\mbox{by\ (\ref{eqn 1.0c})})  \\
&=& m\frac{e^{\mu t}}{1-t}\frac{\partial^{m-1}}{\partial t^{m-1}}\frac{e^{(\lambda-1)t}}{1-t}+(\lambda-1)^m\frac{e^{(\lambda+\mu-1)t}}{1-t}.
\end{eqnarray*}
By (\ref{eqn 2.3}), comparing the coefficient of $\frac{t^n}{n!}$, we get (\ref{eqn 5.4}). By (\ref{eqn 2.4}), $f_{n-k}(\mu+1)$ can be
represented umbrally as $(\mathbf{D}+\mu+1)^{n-k}$, which means that (\ref{eqn 5.4}) is equivalent to
(\ref{eqn 5.3}) by the definition of $Q_{n,m}(\lambda, \mu)$. \qed

Setting $m=2$ in (\ref{eqn 5.4}), by (\ref{eqn 4.3}) and (\ref{eqn 4.5}), we obtain
\begin{theorem}
For any integer $n\geq 0$ and any indeterminants $\lambda, \mu$, there holds
\begin{eqnarray*}
\sum_{k=0}^{n}\binom{n}{k}f_{k+2}(\lambda)\mu^{n-k} &=& \sum_{k=0}^{n}\binom{n}{k}(\lambda^2+2k+1)(\lambda+k)^{k}(\mu-k-1)^{n-k}.
\end{eqnarray*}
\end{theorem}

In general, it seems to be not easy to derive the explicit formula similar to (\ref{eqn 4.3}) and (\ref{eqn 4.4}) for $Q_{n,m}(\lambda, \mu)$,
we leave it as an open problem to the interested readers. One can also be asked to give combinatorial proofs for Corollary 4.2 and 4.6.

%%%%%%%%%%%%%%%%%%%%%%%%%%%%%%%%%%%%%%%%%%%%%%%%%%%%%%%%%%%%%%%%%%%%%%%%%%%%%%%%%%%%%%%%%%%%%%%%%%%%%%%%%%%%%%%%%%%%%%%%
%%% \begin{theorem} For any integers $n,m\geq 0$ and any indeterminants $\lambda, \mu$, there holds                  %%%
%%% \begin{eqnarray}                                                                                                 %%%
%%% \lefteqn{\sum_{k=0}^{n}\binom{n}{k}f_{k+m}(\lambda)(n+1)^{n-k} }\\                                               %%%
%%% &=& (n+1)U_{m-1}(\lambda,n+1)f_{n}(\lambda+n+1)+U_{m}(\lambda,n+1)(n+\lambda)^{n+1},  \label{eqn 4.6}            %%%
%%% \end{eqnarray}                                                                                                   %%%
%%% where $U_{-1}(\lambda, \mu)=\mu^{-1}$, $U_{0}(\lambda, \mu)=0$ and $U_{m}(\lambda, \mu)$ are defined for $m      %%%
%%% \geq 1$ by                                                                                                       %%%
%%% \begin{eqnarray*}                                                                                                %%%
%%% U_{m}(\lambda, \mu)=\sum_{k=0}^{[(m-1)/2]}\binom{m-k-1}{k}\mu^{k}\lambda^{m-2k-1}.                               %%%
%%% \end{eqnarray*}                                                                                                  %%%
%%% \end{theorem}                                                                                                    %%%
%%%%%%%%%%%%%%%%%%%%%%%%%%%%%%%%%%%%%%%%%%%%%%%%%%%%%%%%%%%%%%%%%%%%%%%%%%%%%%%%%%%%%%%%%%%%%%%%%%%%%%%%%%%%%%%%%%%%%%%%

\section*{Acknowledgements} The authors are grateful to the
anonymous referees for the helpful suggestions and comments. The
work was supported by The National Science Foundation of China
(Grant No. 10801020) and supported by the Fundamental
Research Funds for the Central Universities (Grant No. 2009QN070 and 2009QN071).

%==============================================================================================================

\end{document}